\documentclass[10pt,notitlepage]{amsart}
\usepackage{amsmath}     
\usepackage{amsfonts,amsthm}
\usepackage{amssymb}

 \makeatletter 
 \def\activeat#1{\csname @#1\endcsname}
 \def\def@#1{\expandafter\def\csname @#1\endcsname}
 {\catcode`\@=\active \gdef@{\activeat}}

\let\ssize\scriptstyle
\newdimen\ex@   \ex@.2326ex

 \def\requalfill{\cleaders\hbox{$\mkern-2mu\mathord=\mkern-2mu$}\hfill
  \mkern-6mu\mathord=$}
 \def\eqfill{$\m@th\mathord=\mkern-6mu\requalfill}
 \def\deffill{\hbox{$:=$}$\m@th\mkern-6mu\requalfill}
 \def\fiberbox{\hbox{$\vcenter{\hrule\hbox{\vrule\kern1ex
     \vbox{\kern1.2ex}\vrule}\hrule}$}}

 \newdimen\arrwd 
 \newdimen\minCDarrwd \minCDarrwd=2.5pc
        
 \def\findarrwd#1#2#3{\arrwd=#3%
  \setbox\z@\hbox{$\ssize\;{#1}\;\;$}%
  \setbox\@ne\hbox{$\ssize\;{#2}\;\;$}%
  \ifdim\wd\z@>\arrwd \arrwd=\wd\z@\fi
  \ifdim\wd\@ne>\arrwd \arrwd=\wd\@ne\fi}
 \newdimen\arrowsp\arrowsp=0.375em      
 \def\findCDarrwd#1#2{\findarrwd{#1}{#2}{\minCDarrwd}
    \advance\arrwd by 2\arrowsp}
 \newdimen\minarrwd 
 \setbox\z@\hbox{$\longrightarrow$} \minarrwd=\wd\z@

 \def\harrow#1#2#3#4{{\minarrwd=#1\minarrwd
   \findarrwd{#2}{#3}{\minarrwd}\kern\arrowsp
    \mathrel{\mathop{\hbox to\arrwd{#4}}\limits^{#2}_{#3}}\kern\arrowsp}}

 \def@]#1>#2>#3>{\harrow{#1}{#2}{#3}\rightarrowfill}
 \def@>#1>#2>{\harrow1{#1}{#2}\rightarrowfill}
 \def@<#1<#2<{\harrow1{#1}{#2}\leftarrowfill}
 \def@={\harrow1{}{}\eqfill}
 \def@:#1={\harrow1{}{}\deffill}
 \def@ N#1N#2N{\vCDarrow{#1}{#2}\UpDownarrow}
 \def\UpDownarrow{\uparrow\,\Big\downarrow}

\def\hookrightarrowfill{\hbox{$\lhook\joinrel$}\rightarrowfill}
\def@{hk>}#1>#2>{\harrow1{#1}{#2}\hookrightarrowfill}
\def@ c>#1>#2>{\harrow1{#1}{#2}\hookrightarrowfill}
\def\hookleftarrowfill{\leftarrowfill\hbox{$\joinrel\rhook$}}
\def@{hk<}#1<#2<{\harrow1{#1}{#2}\hookleftarrowfill}

 \def@.{\ifodd\row\relax\harrow1{}{}\hfill
   \else\vCDarrow{}{}.\fi}
 \def@|{\vCDarrow{}{}\Vert}
 \def@ V#1V#2V{\vCDarrow{#1}{#2}\downarrow}
 \def@ A#1A#2A{\vCDarrow{#1}{#2}\uparrow}
 \def@(#1){\arrwd=\csname col\the\col\endcsname\relax
   \hbox to 0pt{\hbox to \arrwd{\hss$\vcenter{\hbox{$#1$}}$\hss}\hss}}

 \def\squash#1{\setbox\z@=\hbox{$#1$}\finsm@@sh}
\def\finsm@@sh{\ifnum\row>1\ht\z@\z@\fi \dp\z@\z@ \box\z@}

 \newcount\row \newcount\col \newcount\numcol \newcount\arrspan
 \newdimen\vrtxhalfwd  \newbox\tempbox

 \def\innernewdimen{\alloc@1\dimen\dimendef\insc@unt}
 \def\measureinit{\col=1\vrtxhalfwd=0pt\arrspan=1\arrwd=0pt 
   \setbox\tempbox=\hbox\bgroup$}
 \def\setinit{\col=1\hbox\bgroup$\ifodd\row
   \kern\csname col1\endcsname
   \kern-\csname row\the\row col1\endcsname\fi}
 \def\findvrtxhalfsum{$\egroup
  \expandafter\innernewdimen\csname row\the\row col\the\col\endcsname
  \global\csname row\the\row col\the\col\endcsname=\vrtxhalfwd
  \vrtxhalfwd=0.5\wd\tempbox
  \global\advance\csname row\the\row col\the\col\endcsname by \vrtxhalfwd 
  \advance\arrwd by \csname row\the\row col\the\col\endcsname
  \divide\arrwd by \arrspan
  \loop\ifnum\col>\numcol \numcol=\col%
     \expandafter\innernewdimen \csname col\the\col\endcsname
     \global\csname col\the\col\endcsname=\arrwd
   \else \ifdim\arrwd >\csname col\the\col\endcsname
      \global\csname col\the\col\endcsname=\arrwd\fi\fi
   \advance\arrspan by -1 %
   \ifnum\arrspan>0 \repeat}
 \def\setCDarrow#1#2#3#4{\advance\col by 1 \arrspan=#1 
    \arrwd= -\csname row\the\row col\the\col\endcsname\relax
    \loop\advance\arrwd by \csname col\the\col\endcsname
     \ifnum\arrspan>1 \advance\col by 1 \advance\arrspan by -1%
     \repeat
    \squash{\mathop{
     \hbox to\arrwd{\kern\arrowsp#4\kern\arrowsp}}\limits^{#2}_{#3}}}
 \def\measureCDarrow#1#2#3#4{\findvrtxhalfsum\advance\col by 1%
   \arrspan=#1\findCDarrwd{#2}{#3}%
    \setbox\tempbox=\hbox\bgroup$}
 \def\vCDarrow#1#2#3{\kern\csname col\the\col\endcsname
    \hbox to 0pt{\hss$\vcenter{\llap{$\ssize#1$}}%
     \Big#3\vcenter{\rlap{$\ssize#2$}}$\hss}\advance\col by 1}

 \def\setCD{\def\harrow{\setCDarrow}%
  \def\\{$\egroup\advance\row by 1\setinit}
  \m@th\lineskip3\ex@\lineskiplimit3\ex@ \row=1\setinit}
 \def\endsetCD{$\egroup}
 \def\drop#1\\{\findvrtxhalfsum\advance\row by 2 \measureinit}
 \def\measure{\bgroup
  \def\harrow{\measureCDarrow}%
  \def\\##1{\ifx##1\endmeasure\endmeasure\else\expandafter\drop\fi}%
  \row=1\numcol=0\measureinit}
 \def\endmeasure{\findvrtxhalfsum\egroup}

 \newcount\savedcount   
 \def\LCD#1\end{\savedcount=\count11
   \measure#1\endmeasure
   \vcenter{\setCD#1\endsetCD\kern\medskipamount}%
   \global\count11=\savedcount\end}

\newenvironment{CD}{\let\at=@\catcode`\@=\active\LCD}{\catcode`\@=12\relax}
\makeatother 
\font\smallrm=cmr8

\renewcommand{\c}{\mathcal}
\newcommand{\Og}{\Omega}

\renewcommand{\H}{\text{\rm H}}
\newcommand{\h}{\text{\rm h}}
\newcommand{\Hom}{\text{\rm Hom}}
\newcommand{\IHom}{\it Hom}

\newcommand{\codim}{\text{\rm codim}}
\newcommand{\wt}{\widetilde}
\renewcommand{\:}{\colon}
\newcommand{\x}{\times}
\newcommand{\ox}{\otimes}
\newcommand{\bx}{\boxtimes}
\newcommand{\CIP}{\text{\bf P}_{\text{\bf C}}}
\newcommand{\IP}{\text{\bf P}}

\renewcommand{\deg}{\mathrm{deg}\,}
\renewcommand{\dim}{\mathrm{dim}\,}
\newcommand{\reg}{\mathrm{reg}\,}
\newcommand{\smashedrightarrow}
{\setbox0=\hbox{$\rightarrow$}\ht0=1pt\box0}
\newcommand{\risom}{\buildrel{\hskip-0.08cm\sim}\over{\smashedrightarrow}}

\newcommand{\bigomega}{\hbox{\large $\omega$}}

\let\bog=\w
\def~{\thinspace\nobreak}

\newtheorem{theorem}{Theorem}[section]
\newtheorem{corollary}[theorem]{Corollary}
\newtheorem{lemma}[theorem]{Lemma}
\newtheorem{proposition}[theorem]{Proposition}
\theoremstyle{definition}
\newtheorem{definition}[theorem]{Definition}
\newtheorem{remark}[theorem]{Remark}
\newtheorem{subsct}[theorem]{}
\theoremstyle{plain}

\begin{document}
\title[{\smallrm Bounding solutions of Pfaff equations}] {Bounding
  solutions of Pfaff equations}
 \author[{\smallrm E. Esteves and S. Kleiman}]{E. Esteves
 \ and \ S. Kleiman} 

\date{}
\maketitle

\begin{center}
Instituto de Matem\'atica Pura e Aplicada\\
Estrada Dona Castorina 110,\\
 22460-320 Rio de Janeiro RJ, Brazil\\
E-mail: \texttt{esteves@impa.br}\\
\medskip

Department of Mathematics, MIT \\
77 Massachusetts Avenue \\
Cambridge, MA 02139, USA \\
E-mail: \texttt{kleiman@math.mit.edu}
\end{center}
\bigskip

\begin{abstract} 
 Let $\omega$ be a Pfaff system of differential forms on $\CIP^n$.  Let
$S$ be its singular locus, and $Y$ a solution of $\omega=0$.  We prove
$Y\cap S$ is of codimension at most 1 in $Y$, just as Jouanolou
suspected; he proved this result assuming $\omega$ is completely
integrable, and asked if the integrability is, in fact, needed.
Furthermore, we prove a lower bound on the Castelnuovo--Mumford
regularity of $Y\cap S$.  As in two related articles, we derive upper
bounds on numerical invariants of $Y$, thus contributing to the solution
of the Poincar\'e problem.  We work with Pfaff fields not necessarily
induced by Pfaff systems, with ambient spaces more general than
$\CIP^n$, and usually in arbitrary characteristic.  \end{abstract}

\section{Introduction}
 In his seminal work \cite{J} on algebraic Pfaff equations, Jouanolou
proved that a complex analytic foliation of positive dimension of an
open subset $U\subseteq\CIP^n$ has no compact leaves; see his Prop.~4.2,
p.~130.  He went on, in his Cor.~4.2.7, p.~133, to prove that, if the
foliation arises from a {\it completely integrable\/} system of Pfaff
forms, if $U$ is the complement of the singular locus $S$ of the system,
and if $Y\subseteq\CIP^n$ is a closed subvariety such that $Y\cap U$ is
a leaf, then $S$ intersects $Y$ in codimension at most 1.  Then in
Rem.~4.2.8, p.~134, he said it would be interesting to remove the
hypothesis of complete integrability.
 \renewcommand{\thefootnote}{}%
 \footnote 
  {2000 {\it Mathematics Subject Classification} 37F75 (primary), 
32S65, 14H99, 14B05, 14F10 (secondary).} 
 \footnote 
  {{\it Keywords} Foliations, Pfaff equations, Poincar\'e problem,
singularities, C-M regularity.}

The present article advances Jouanolou's work.  As he envisioned, it is
indeed possible to work with an arbitrary Pfaff system $\omega$: without
assuming integrability, we prove that the singular locus $S$ of $\omega$
intersects a solution $Y$ of $\omega=0$ in a subvariety of codimension
at most 1.

In fact, we go further.  Let $X$ be, more generally, a complex
projective scheme of pure dimension $n$.  A Pfaff system on $X$ induces,
via exterior powers and the perfect pairing of differential forms, a map
$\eta\:\Og^b_X\to\c L$ from the sheaf of differential forms to an
invertible sheaf; see Subsec.~3.1.  However, the converse is not
true: such an $\eta$ does not come necessarily from a Pfaff system.

Let $\eta\:\Og^b_X\to\c L$ be a nonzero map with $0<b<n$.  Its {\it
singular locus\/} is defined as the subscheme $S\subset X$ of points
where $\eta$ is not surjective. Let $Y\subset X$ be a reduced closed
subscheme of dimension $b$.  Assume no $b$-dimensional component lies
in $S$.  Assume $Y$ is {\it invariant\/} under $\eta$; that is, $\eta|Y$
factors through the natural map $\Og^b_X|Y\to\Og^b_Y$.  If $X$ is smooth,
and $\eta$ arises from
a Pfaff system, then $Y$ is a solution of the system in Jouanolou's
sense; see Subsec.~3.1 and Prop.~3.2.

Under the above conditions,  Prop.~3.3 says that the sheaf of ideals
$\c J$ of $Y\cap S$ in $X$ satisfies $\H^b(\c J\ox\c L)\neq 0$ and that,
if
the induced map $\H^b(\eta)\:\H^b(\Og^b_X)\to\H^b(\c L)$ vanishes, then
$Y\cap S$ has codimension 1 in $Y$ and $\h^b(\c L|Y)<\h^b(\Og^b_Y)$.

If $X=\CIP^n$, then $\H^b(\eta)=0$ because $\H^b(\c L)=0$; so
$\codim(Y\cap S,Y)=1$. Furthermore, if $\H^b(\c J\ox\c L)\neq 0$, then
the Castelnuovo--Mumford regularity $\reg(Y\cap S)$ is greater than
$m:=\deg\c L+b$; see Cor.~4.5.

As Soares observed in the introduction to \cite{S}, Jouanolou's work can
be used to tackle the Poincar\'e problem.  Soares' observation served to
motivate our work here and in \cite{EK1} and \cite{EK2}.

In 1891, Poincar\'e \cite{P}, p.~161, posed the problem of bounding the
degree of an algebraic curve $Y$ invariant under a polynomial vector
field on the complex plane.  Versions of this problem have been
considered in a number of recent works; references are given in
\cite{EK1}.  {}From our point of view, the general problem is simply to
find upper bounds on the various numerical invariants of  $Y$.

Roughly, Soares' idea is this: upper bounds on the numerical invariants
of $Y$ arise from lower bounds on the numerical invariants of $Y\cap S$,
where, as always, $S$ is the singular locus of the field.  In
\cite{EK1}, this idea is used to improve bounds obtained by Campillo,
Carnicer, and Garc\'\i a de la Fuente \cite{CCG}, and by Du Plessis and
Wall \cite{dPW}.

In the present article, we approach the Poincar\'e problem in a new way.
It is based on the inequality $\h^b(\c L|Y)<\h^b(\Og^b_Y)$, which
obtains if $\H^b(\eta)=0$ according to Prop.~3.3.  Our Cor.~4.5 gives
one application: if $X=\CIP^n$, if $Y$ is arithmetically
Cohen--Macaulay~---~for instance, a complete intersection~---~and if
$\h^b(\Og^b_Y)=1$, then $\reg(Y)\le m+1$. The third condition
$\h^b(Y,\Og^b_Y)=1$ is satisfied when $Y$ is integral and has
normal-crossings in codimension 1; see Rem.~4.7.  Since the regularity
of a plane curve is just its degree, we recover a fundamental result
proved by Cerveau and Lins Neto \cite{CL}.  We also recover
\cite{E}, Thm.~1, p.~3, which generalizes their result to curves in
$\CIP^n$.

The condition $\h^b(Y,\Og^b_Y)=1$ is also satisfied when $Y$ has higher
singularities, yielding new solutions to the original Poincar\'e problem
on $\CIP^2$.  Indeed, assume $Y$ is a plane curve of degree $d$. Let
$\Sigma$ be its singular locus, the subscheme cut out by its
polars.  Set $\sigma:=\reg(\Sigma)$.  In \cite{EK2}, Thm.~2.5 asserts
that $d\leq m+1$ if $\sigma\leq d-2$; otherwise, $2d\leq m+\sigma+3$,
with equality if $d\geq 2m$ and $S$ is finite.

Our Prop.~3.3 applies to ambient varieties other than projective space.
For instance, it applies to multiprojective space; see Thm.~4.3.

Proceeding in a different direction, assume $X$ is smooth and $\text{\rm
Pic}(X)=\text{\bf Z}$.  We obtain two results. First, Prop.~3.4 says
that, if the normal sheaf of $Y$ in $X$ has positive degree on some
curve lying in the smooth locus of $Y$, then $\codim(Y\cap S,\,Y)=1$.
Second, Thm.~3.6 says that, if $Y$ is a hypersurface with
normal-crossings in codimension 1, then $\deg Y\le\deg\c L(-K)$ where
$K$ is a canonical divisor of $X$.  This theorem generalizes part of the
main theorem in \cite{BM}, p.~594.

Using methods similar to Jouanolou's, Lehmann \cite{Le} too advanced his
work. However, our results seem to be completely independent of
Lehmann's; and our methods, completely different.

Surprisingly, our results rest on a rather unsurprising fact: the map
$\H^b(\Og^b_X)\to\H^b(\Og^b_Y)$ does not vanish.  This nonvanishing was
known in some generality, at a minimum when $X$ and $Y$ are smooth; and
probably it was expected in the generality we need.  However, there
appears to be no suitable reference.  Some references are too abstract;
others, not general enough.  So the fact is proved in Prop.~2.1.

All our schemes are defined over a fixed algebraically closed field.
All our results hold over any field of characteristic 0, not just
$\text{\bf C}$.  Except for Prop.~3.4 and Thm.~3.6, all our results hold
over a field of characteristic $p>0$ if the restriction map
$\H^b(\Og^b_X)\to\H^b(\Og^b_Y)$ does not vanish.  Prop.~2.2 gives
sufficient conditions for this nonvanishing.  For instance, if
$X=\IP^n$, then it is enough that $p\nmid\deg Y$.  The proof of
Prop.~2.2 is similar to that of Prop.~2.1, but is more involved, most
notably in its use of the theory of residues.  We feel the effort is
worthwhile, owing to the resurgence of interest in foliations in
positive characteristic, caused by McQuillan's proof in \cite{M}  of the
Green--Griffiths conjecture, which uses Miyaoka's results proved by
means of reduction to positive characteristic.

\section{Nonvanishing}

 \begin{proposition} Let $X$ be a projective scheme over a field of
characteristic zero, and $f\:Y\to X$ a finite map.  Set $b:=\dim Y$.
Then the natural map $\H^b(\Og^b_X)\to\H^b(\Og^b_Y)$ is nonzero.
 \end{proposition}
 \begin{proof} 
We proceed by induction on $b$.  If $b=0$, then the map in question is
just the pullback map $\H^0(\c O_X)\to\H^0(\c O_Y)$, which is always
nonzero.  So assume $b>0$.

Let $Y'\subseteq Y$ be an irreducible component of dimension $b$, and give
$Y'$  the reduced structure.  It is enough to show the composition
        $$\H^b(\Og^b_X)\to\H^b(\Og^b_Y)\to\H^b(\Og^b_{Y'})$$
 is nonzero.  So we may replace $Y$ by $Y'$, and thus assume $Y$ is
integral.

Let $\pi\: Y^*\to Y$ be the normalization map. It is enough to show the 
natural map $\H^b(\Og^b_X)\to\H^b(\Og^b_{Y^*})$ is nonzero. Since $\pi$ is 
finite and, hence, $\dim Y^*=b$, we may replace $Y$ by
$Y^*$, and thus assume $Y$ is normal.

Let us now find on $X$ an effective Cartier divisor $E$ satisfying the
following conditions:\par
 \smallskip
 (2.1.1)\enspace The preimage $F:=f^{-1}(E)$ is nonempty, Cartier and
smooth in
 \indent\hphantom{(4.1.3)}\enspace codimension~1.\par
 (2.1.2)\enspace No component of the singular locus of $Y$ is contained
in $F$.\par
 (2.1.3)\enspace  The induced maps $\Og^b_X\to\Og^b_X(E)$ and
$\Og^b_Y\to\Og^b_Y(F)$ are injective.\par
  (2.1.4)\enspace If $b>1$, then  $\H^{b-1}(\Og^b_Y(F))=0$.
 \smallskip

To start, let $E$ be any effective very ample divisor such that
$E\not\supset f(Y)$.  However, if $b>1$, then take $E$ ample enough so
that $\H^{b-1}(f_*\Og^b_Y(E))=0$.  Then $F:=f^{-1}(E)$ is Cartier on
$Y$.  Hence (2.1.4) holds.  Moreover, $F$ is nonempty because $E$ is
ample and $\dim f(Y)=b>0$.

Vary $E$ inside its complete linear system, keeping $E\not\supset f(Y)$.
Correspondingly, $F$ traces on $Y$ a linear system without base points
(although it may be incomplete).  If $E$ is general, then $E$ and $F$
contain no associated point of $\Og^b_X$ and $\Og^b_Y$ respectively;
hence, (2.1.3) holds. Similarly, (2.1.2) holds if $E$ is general.

Finally, since the characteristic is 0, if $E$ is general,
then $F$ is smooth off the singular locus of $Y$ by a form of Bertini's 
first theorem; see \cite{K}, Cor.~5, p.~291.  In particular, 
$F$ is smooth in codimension 1 by (2.1.2). Then (2.1.1) holds.

 Consider now the second fundamental exact sequence:
        $$\c O_X(-E)|E\to\Og^1_X|E\to\Og^1_E\to 0.$$ 
 In a standard way, it induces a map,
        $$\eta_{E,X}\:\Og^{b-1}_E\to\Og^b_X(E)|E;$$
 namely, given the germ of a form on $E$, lift it to $X$, then
wedge with the meromorphic 1-form $dt/t$ where $t=0$ is a local equation
for $E$, and finally restrict to $E$.

Tensoring the standard exact sequence
        $$0\to\c O_X \to\c O_X(E)\to\c O_X(E)|E\to 0$$
 with $\Og^b_X$, we obtain a sequence
        $$0\to\Og^b_X\to\Og^b_X(E)\to\Og^b_X(E)|E\to 0,$$
 which is exact on the left by (2.1.3) above.  Form the coboundary map
        $$h_{E,X}\:\H^{b-1}(\Og^b_X(E)|E)\to\H^b(\Og^b_X),$$
 and set $v_{E,X}:=h_{E,X}\circ\H^{b-1}(\eta_{E,X})$.

Similarly, for $Y$ and $F$, we have maps $\eta_{F,Y}$, $h_{F,Y}$ and
$v_{F,Y}$.  Form the diagram
  $$
  \begin{CD}
  \H^{b-1}(\Og^{b-1}_E) @>>>\H^{b-1}(\Og^{b-1}_F) \\
  @VVv_{E,X}V @VVv_{F,Y}V\\
  \H^b(\Og^b_X) @>>> \H^b(\Og^b_Y)
  \end{CD}
  \eqno(2.1.5)
  $$
 using the natural horizontal maps.  It is plainly commutative.

By induction, the top map is nonzero.  
Now, $F$ is smooth in codimension 1 by (2.1.1).  
In addition, the singular locus of $Y$ intersects 
$F$ in codimension 2 by (2.1.2). Therefore,
$\eta_{F,Y}\:\Og^{b-1}_F\to\Og^b_Y(F)|F$ is an isomorphism in
codimension 1.  Hence $\H^{b-1}(\eta_{F,Y})$ is an isomorphism.

First assume $b>1$.  Then $\H^{b-1}(\Og^b_Y(F))=0$ by (2.1.4). So
the coboundary map
        $$h_{F,Y}\:\H^{b-1}(\Og^b_Y(F)|F)\to\H^b(\Og^b_Y)$$
 is injective.  Hence, in Diagram (2.1.5), the top-right composition is
nonzero.  Hence the left-bottom composition is also.  Therefore, the
bottom map is nonzero.

Finally, assume  $b=1$.  In this case,  Diagram (2.1.5) becomes
  $$
  \begin{CD}
  \H^0(\c O_E) @>>> \H^0(\c O_F)\\
  @VVv_{E,X}V @VVv_{F,Y}V\\
  \H^1(\Og^1_X) @>>> \H^1(\Og^1_Y)
  \end{CD}
  \eqno(2.1.6)
  $$
 As before, we need only show that the top-right composition is nonzero.
To do so, we need only prove $v_{F,Y}(1)\neq 0$.

By definition, $v_{F,Y}$ is the following composition:
  $$
  \begin{CD}
  \H^0(\c O_F) @>\H^0(\eta_{F,Y})>> \H^0(\Og^1_Y(F)|F) @>h_{F,Y}>> 
  \H^1(\Og^1_Y).
  \end{CD}
  $$
 Given $y\in F$, let $t$ be a uniformizing parameter of $F$ at $y$.  Then
$\eta_{F,Y}(1)$ is at $y$ equal to the class of $dt/t$.  Now, let
$\rho_Y\:\H^1(\Og^1_Y)\to k$ be the global residue map; we compute it
by summing local residues.  Consequently, $\rho_Y(v_{F,Y}(1))=\deg F$.
Since $F$ is nonempty, $\deg F\neq0$.  Hence $v_{F,Y}(1)\neq 0$, and the
proof is complete.
 \end{proof}

 \begin{proposition} Let $X$ be a projective scheme over a field of
characteristic $p>0$, and $f\:Y\to X$ a finite map.  Set $b:=\dim Y$.
Assume there are Cartier divisors $E_1,\dots,E_b$ on $X$ such that
        $$\int_X E_1\cdots E_b\cdot f_*[Y]\not\equiv 0\text{ \rm (mod
$p$)}.
        \eqno(2.2.1)$$
 Then the natural map $\H^b(\Og^b_X)\to\H^b(\Og^b_Y)$ is nonzero.  
 \end{proposition}

\begin{proof}
The proof is analogous to that of Proposition 2.1 (and reproves the 
proposition); we dwell only on the alterations.  They are required 
because we can no longer guarantee $F$ is smooth in codimension 1.  
Notably, we must use more of the theory of residues.

As before, we may assume that $b>0$.  Again, we may replace $Y$ by some 
integral component $Y'$; indeed, (2.2.1) will still hold as $[Y]$ is a 
linear combination of the $[Y']$ of dimension $b$. Then $f$ is 
generically \'etale; indeed, if $n:=\deg f$, then $f_*[Y]=n[f(Y)]$, and so 
$p\nmid n$ owing to (2.2.1).

So we may assume that $Y$ is generically smooth of pure dimension $b$
and that $f$ is generically \'etale onto its image.  We are going to
prove a stronger assertion, namely the nonvanishing of the composition
$$
  \begin{CD}
        \H^b(\Og^b_X)@>>>\H^b(\Og^b_Y)@>\rho_Y>>k.$$
  \end{CD}
  \eqno(2.2.2)
$$
Here $\rho_Y$ is the generalized residue map, defined as explained in the 
next paragraph.

Given an integral, projective scheme $Z$ of dimension $e$, 
let us denote by $\rho_Z\:\H^e(\Og^e_Z)\to k$ its generalized residue map; 
see Thm.~0.1 on p.~10 of \cite{Li} and the discussion thereafter, where 
$\rho_Z$ is denoted by $\int_Z$ however. Given a generically smooth, 
projective scheme $Z$ of pure dimension $e$, let 
$\rho_Z$  denote the composition
$$
\begin{CD}
\H^e(\Og^e_Z)\to\H^e(\Og^e_{Z_1})\oplus\cdots\oplus\H^e(\Og^e_{Z_s})
@>(\rho_{Z_1},\dots,\rho_{Z_s})>> k,
\end{CD}
$$
where $Z_1,\dots,Z_s$ are the irreducible components of $Z$ with their 
reduced structures, and the first map is the natural one.

As in the proof of Proposition 2.1, we can find an effective very
ample divisor $E$ on $X$ such that $F:=f^{-1}(E)$ is Cartier, nonempty, 
and (2.1.2) and (2.1.3) hold.  In addition, as we are going to see, we may 
assume the following three conditions hold:\par
 \smallskip
 (2.2.3)\enspace We have $\int E\cdot E_2\cdots E_b\cdot f_*[Y] \not\equiv
0\text{ \rm (mod $p$)}$.\par
 (2.2.4)\enspace The scheme $F$ is generically smooth and $f|F$ is generically 
\'etale
 \indent\hphantom{(4.1.3)}\enspace onto its image.\par
 (2.2.5)\enspace There are a finite map $g\:Y\to P$, where $P:={\bf
  P}^b$, and a hy-
 \indent\hphantom{(4.1.3)}\enspace perplane $M\subset P$ such that
 $g^{-1}M=F$ and $g|F$ is generically
 \indent\hphantom{(4.1.3)}\enspace \'etale onto $M$.
 \smallbreak

If (2.2.3) doesn't already hold, then replace $E$ by a general member
of the linear system $|mE+E_1|$ for $m\gg0$.  Then (2.1.2), (2.1.3) and
(2.2.3) hold.

As to (2.2.4), since generically $f$ is \'etale and $Y$ is smooth, 
$Y$ has a smooth, dense open subset $U$ such that $f|U$ is \'etale over 
$f(Y)$. 
We may replace $E$ by a general member of the linear system $|E|$, and 
assume that every component of $F$ intersects $U$. Furthermore, even 
though $p>0$, we may assume $F\cap U$ is smooth by another form of
Bertini's first theorem; see \cite{K}, Cor.~12, p.~296. Then $F$ is 
generically smooth and $f|F$ is generically \'etale onto $Z:=f(F)$.

To ensure (2.2.5),  use the system $|E|$ to embed $X$ in a
projective space $P'$, and let $E'$ be a hyperplane that cuts $E$ out of
$X$.  Let $z_1,\dots,z_s\in Z$ be simple points, one for 
each component of $Z$.  Let $C'$ be a linear subspace of $E'$ of 
codimension $d$ 
such that $C'$ misses both $Z$ and its tangent spaces 
$T_{z_i}Z\subset E'$.  Then $C'$ misses $f(Y)$ too.  So projection from
$C'$ 
induces a finite map $g'\:f(Y)\to{\bf P}^b$. Set $P:={\bf P}^b$ and 
$g:=g'f$.

There is a hyperplane $M\subset P$ such that ${g'}^{-1}M=Z$ since
$C'\subset E'$.  Hence $g'|Z$ is finite onto $M$.  It is also \'etale 
at each $z_i$ since $C'$ misses $T_{z_i}Z$. In particular, 
$g'|Z$ is generically \'etale onto $M$. Since also $f|F$ is generically 
\'etale onto its image, by (2.2.4), the composition $g|F$ is 
generically \'etale onto $M$. Thus (2.2.5) holds.

We proceed by induction on $b\ge1$, using the diagrams of maps (2.1.5) and 
(2.1.6), which exist and are commutative by (2.1.3).  

First assume $b=1$.  Since Diagram (2.1.6) is commutative, we need only
prove that $\rho_Y(v_{F,Y}(1))\neq 0$.  Now, $Y$ is smooth along $F$ by
(2.1.2).  Given $y\in F$, let $t$ be a uniformizing parameter of $F$ at
$y$.  Since $F$ is of pure dimension 0 and generically smooth by
(2.2.4), $F$ is reduced.  Then, as before, $\eta_{F,Y}(1)$ is at $y$
equal to the class of $dt/t$.  Consequently, $\rho_Y(v_{F,Y}(1)) = \deg
F$ in $k$.  However, $\deg F = \int E\cdot f_*[Y]$.  Hence (2.2.3)
implies $\rho_Y(v_{F,Y}(1))\neq 0$, as desired.

Finally, assume $b>1$.  Since $F$ is generically smooth with pure
dimension 
$b-1$, and $f|F$ is generically \'etale onto its image, 
by induction the composition
  $$
  \begin{CD}
        \H^{b-1}(\Og^{b-1}_E)@>>>\H^{b-1}(\Og^{b-1}_F)@>\rho_F>>k
  \end{CD}
  $$
 is nonzero.  Now, Diagram (2.1.5) is commutative.  It will follow that
the composition (2.2.2) is nonzero once we prove that the following
diagram is commutative:
  $$
  \begin{CD}
  \H^{b-1}(\Og^{b-1}_F) @>\rho_F>> k\\
   @VVv_{F,Y}V @|\\
  \H^b(\Og^b_Y) @>\rho_Y>> k
  \end{CD}
  $$
 We are going to reduce the matter to the case where $M$ and $P$ replace
$F$ and $Y$.

Using the natural maps, form the following diagram:
  $$
  \begin{CD}
  \H^{b-1}(g_*\c O_F\otimes\Og^{b-1}_M) @>>>
    \H^{b-1}(g_*\Og^{b-1}_F)@=\H^{b-1}(\Og^{b-1}_F) @>\rho_F>> k\\
   @VVV @VVV @VVv_{F,Y}V @|\\
  \H^b(g_*{\c O_Y}\otimes\Og^b_P) @>>>
    \H^b(g_*\Og^b_Y) @= \H^b(\Og^b_Y) @>\rho_Y>> k
  \end{CD}\eqno(2.2.6)
  $$
 A look at the construction of the left-hand square shows it is
commutative.  Its top map is surjective; indeed, $g|F$ is 
generically \'etale onto $M$ by (2.2.5), so 
$g_*\c O_F\otimes\Og^{b-1}_M\to g_*\Og^{b-1}_F$ is
generically surjective.  It will follow that the right-hand square is
commutative once we prove that the outer ``square'' is commutative.

Let $Y_1,\dots,Y_t$ be the irreducible components of $Y$ with their
reduced 
structure. The bottom composition in (2.2.6) is equal to the following 
composition of natural maps:
\[
  \H^b(g_*\c O_Y\otimes\Og^b_P) \xrightarrow{}
  \bigoplus_{i=1}^t\H^b(g_*\c O_{Y_i}\otimes\Og^b_P) \xrightarrow{}
  \bigoplus_{i=1}^t\H^b(\Og^b_{Y_i})
\xrightarrow{(\rho_{Y_1},\dots,\rho_{Y_t})} k.
\]
 By Thm.~0.1(b) on p.~10 of \cite{Li}, for each $i$ the diagram below 
commutes:
$$
\begin{CD}
\H^b(g_*\c O_{Y_i}\otimes\Og^b_P) @>>> \H^b(\Og^b_{Y_i})\\
@VVV @V\rho_{Y_i}VV\\
\H^b(\Og^b_P) @>\rho_P>> k,
\end{CD}
$$
where the left vertical map is induced by the trace map 
$g_*\c O_{Y_i}\to\c O_P$. Since $Y$ is generically smooth, 
the trace map $g_*\c O_Y\to\c O_P$ is the 
sum of the trace maps $g_*\c O_{Y_i}\to\c O_P$. 
So the bottom composition in (2.2.6) is 
equal to the bottom composition below:
$$
  \begin{CD}
  \H^{b-1}(g_*\c O_F\otimes\Og^{b-1}_M) @>>>
    \H^{b-1}(\Og^{b-1}_M) @>\rho_M>> k\\
   @VVV @VVv_{M,P}V @|\\
  \H^b(g_*{\c O_Y}\otimes\Og^b_P) @>>> \H^b(\Og^b_P) @>\rho_P>> k
  \end{CD}\eqno(2.2.7)
  $$
 where the left-hand horizontal maps are induced by the two trace maps
$g_*\c O_F\to\c O_M$ and $g_*\c O_Y\to\c O_P$.  The latter map restricts
to the former, and it follows that the left-hand square is commutative.
By analogy, the top composition in (2.2.6) is equal to that in 
(2.2.7). So the outer ``square'' in (2.2.6) is commutative if 
the right-hand square in (2.2.7) is commutative.

By Thm.~0.1(a) on p.~10 of \cite{Li}, $\rho_M$ and $\rho_P$ are the
``well-known canonical isomorphisms.''  A simple explicit calculation
now shows the right-hand square in (2.2.7) is commutative.  
\end{proof}

\section{Pfaff fields}

\begin{subsct}
 \emph{Pfaff systems, equations, and fields.}  A \emph{Pfaff system of
rank} $a$ on a smooth scheme $X$ of pure dimension $n$ over a field is,
according to Jouanolou \cite{J}, pp.~136--38, a nonzero map $u\:\c
E\to\Og^1_X$ where $\c E$ is a locally free sheaf of constant rank $a$
with $0<a<n$.  The \emph{singular locus} of the system $u$ is the closed
subscheme $S$ of $X$ whose ideal $\c I_S$ is the image of the induced
map $\bigwedge^a\c E\ox(\Og^a_X)^*\to\c O_X$.  A \emph{solution} is a
closed subscheme $Y$ of $X$ with pure codimension $a$ such that the map
 $$\textstyle
 (\bigwedge^a d)\wedge(\bigwedge^a u|Y)\:
 \bigwedge^a(\c I_{Y,X}/\c I_{Y,X}^2)\ox\bigwedge^a\c E|Y
 \to\bigwedge^2\Og^a_X|Y
 $$
 vanishes, where $d\:\c I_{Y,X}/\c I_{Y,X}^2\to\Og^1_X|Y$ is the standard 
map, given by differentiation.

The notions of singular locus and solution involve the map $\bigwedge^a
u$, not $u$ directly.  So it is natural to generalize the theory in the
following way; compare with Brunella and Mendes \cite{BM}, pp.~593--94.
Define a \emph{Pfaff equation of rank} $a$ to be an equation $\sigma=0$
where $\sigma$ is a nonzero global section of $\Og^a_X\otimes\c
N$ for a given integer $a$ and a given invertible sheaf $\c N$.  The
\emph{singular locus} is the closed subscheme whose ideal is the image
of the dual map $(\Og^a_X)^*\otimes\c N^*\to\c O_X$.  A \emph{solution}
is a closed subscheme $Y$ with pure codimension $a$ such that the
following
natural map vanishes:
 $$\textstyle
 (\bigwedge^a d)\wedge(\sigma\ox\c N^*)|Y\:
 \bigwedge^a(\c I_{Y,X}/\c I_{Y,X}^2)\ox\c N^*|Y\to\bigwedge^2\Og^a_X|Y.
 \eqno(3.1.1)
 $$

Alternatively, we may view a Pfaff equation as follows.  Set $b:=n-a$.
Let $\tau\:\Og^a_X\ox\Og^b_X\to\Og^n_X$ be the natural pairing; $\tau$
is perfect since $X$ is smooth of pure dimension $n$.  Set $\c
L:=\Og^n_X\ox\c N$.  Then there is a natural isomorphism,
       $$\H^0(\Og^a_X\ox\c N)=\Hom(\Og^b_X,\c L),$$
  under which $\sigma$ corresponds to the composition
        $\eta:=(\tau\ox\c N)\circ(\sigma\ox\Og^b_X)$. 
 Thus giving $(a,\,\c N,\,\sigma)$ is equivalent to giving $(b,\, \c
L,\,\c \eta)$.
 
More generally, without assuming that $X$ is smooth or equidimensional,
define a \emph{Pfaff field of rank} $b$ to be a nonzero map
$\eta\:\Og^b_X\to\c L$ for a given integer $b$ with $0<b<n$ and a given
invertible sheaf $\c L$.  Define the \emph{singular locus} $S$ of $\eta$
to be the closed subscheme of $X$ whose ideal $\c I_S$ is the image of
the induced map $\Og^b_X\ox\c L^*\to\c O_X$.  Say that a closed
subscheme $Y$ of $X$ is \emph{invariant} under $\eta$ if the restriction
$\eta|Y\:\Og^b_X|Y\to\c L|Y$ factors through the standard map
$\beta\:\Og^b_X|Y\to\Og^b_Y$, in other words, if there is a commutative
diagram
  $$
  \begin{CD}
  \Og^b_X @>\eta>> \c L\\
  @V VV @VVV\\
  \Og^b_Y @>\mu>> \c L|Y
  \end{CD}
  \eqno(3.1.2)
  $$
 whose vertical maps are the standard maps.

Again assume that $X$ is smooth of pure dimension $n$.  Then it is easy
to see that the singular locus of the Pfaff field $\eta$ is the same as
the singular locus of the corresponding Pfaff equation $\sigma=0$.  Now,
to avoid uninteresting cases, assume $Y$ is reduced and equidimensional
and no component lies in $S$.  Then $Y$ is invariant under $\eta$ if and
only if $Y$ is a solution of $\sigma=0$, at least when $\dim Y=b$,
according to the next proposition.
\end{subsct}

 \begin{proposition} Let $X$ be a smooth equidimensional scheme,
$\eta\:\Og^b_X\to\c L$ a Pfaff field, $S$ its singular locus, and
$\sigma=0$ the corresponding Pfaff equation.  Let $Y$ be a closed
subscheme; assume $Y$ is reduced of pure dimension $b$ and no component
lies in $S$.  Then $Y$ is invariant under $\eta$ if and only if $Y$ is a
solution of $\sigma=0$.
 \end{proposition}

 \begin{proof} Since no component of $Y$ lies in
$S$, there is an open subset $U$ of $X-S$ such that $V:=U\cap Y$ is
dense in $Y$.  Since $Y$ is reduced, we may assume $V$ is smooth.
Moreover, $Y$ is invariant under $\eta$ if and only if $V$ is invariant
under $\eta|U$; indeed, if we let $\c K$ denote the kernel of the map
$\Og^b_X|Y\to\Og^b_Y$, then $(\eta|Y)(\c K)$ vanishes if and only if its
restriction over $V$ vanishes, since $\c L$ is invertible and $Y$ is
reduced.  Similarly, $Y$ is a solution of $\sigma=0$ if and only if $V$
is a solution of $\sigma|U=0$; indeed, the image of the map in (3.1.1)
vanishes if and only if its restriction over $V$ vanishes.  Replacing
$X$ by $U$ and $Y$ by $V$, we may thus assume $Y$ is smooth and $S$ is
empty.

Let $n:=\dim X$ and $a:=n-b$. Set $\c N:=\c L\ox(\Og^n_X)^{-1}$ and 
$\c J:=\c I_{Y,X}$. 
Since $X$ and $Y$ are smooth, $d\:\c J/\c J^2\to\Og^1_X|Y$ is
locally split injective; whence, so is $\bigwedge^a d$.  Since $S$ is
empty, $\sigma\ox\c N^*\:\c N^*\to\Og^a_X$ is locally split injective
too; whence, so is its restriction to $Y$.  It follows that the map in
(3.1.1) vanishes if and only if there is a map, necessarily an
isomorphism, $\zeta\:\c N^*|Y\to\bigwedge^a(\c J/\c J^2)$ such that
$\bigwedge^a d\circ\zeta$ is equal to $(\sigma\ox\c N^*)|Y$.  In other
words, $Y$ is a solution of $\sigma=0$ if and only if such a $\zeta$
exists.

On the other hand, consider the map $\beta\:\Og^b_X|Y\to\Og^b_Y$, and
form the composition
 $$
 \begin{CD}
 \alpha\:\c N^*\ox\Og^b_X @>(\sigma\ox\c N^*)\ox\Og^b_X>> 
 \Og^a_X\ox\Og^b_X @>\tau>> \Og^n_X
 \end{CD}
 $$
 where $\tau$ is the natural pairing of forms.  Virtually by definition,
$Y$ is invariant under $\eta$ if and only if there is a map $\gamma\:\c
N^*|Y\ox\Og^b_Y\to\Og^n_X|Y$ such that $\alpha|Y=\gamma\circ(\c
N^*|Y\ox\beta)$.

Since $X$ and $Y$ are smooth, the natural map
 $ \psi\:\bigwedge^a(\c J/\c J^2)\ox\Og^b_Y\to\Og^n_X|Y$
  is an isomorphism such that
 $$
 \textstyle
 \tau|Y\circ\bigl(\bigwedge^a d\ox\Og^b_X|Y\bigr)=
 \psi\circ\bigl(\bigwedge^a(\c J/\c J^2)\ox\beta\bigr).\eqno(3.2.1)
 $$
 Since $\Og^b_Y$ is invertible, every map $\gamma\:\c
N^*|Y\ox\Og^b_Y\to\Og^n_X|Y$ can be written in the form
$\gamma=\psi\circ(\zeta\ox\Og^b_Y)$ where $\zeta\:\c
N^*|Y\to\bigwedge^a(\c J/\c J^2)$.  Hence $Y$ is invariant under $\eta$
if and only if there exists a $\zeta$ such that
$\alpha|Y=\psi\circ(\zeta\ox\Og^b_Y)\circ(\c N^*|Y\ox\beta)$.

Given $\zeta$, using the funtoriality of $\ox$ twice and Equation
(3.2.1), we obtain
 \begin{align*}
 \psi\circ(\zeta\ox\Og^b_Y)\circ(\c N^*|Y\ox\beta)
 &=\textstyle\psi\circ\bigl(\bigwedge^a(\c J/\c J^2)\ox\beta\bigr)
        \circ(\zeta\ox\Og^b_X|Y)\\
 &=\textstyle\tau|Y\circ\bigl(\bigwedge^a d\ox\Og^b_X|Y\bigr)
        \circ(\zeta\ox\Og^b_X|Y)\\
 &=\textstyle\tau|Y\circ(\bigl(\bigwedge^a d\circ\zeta\bigr)
        \ox\Og^b_X|Y).
 \end{align*}
 Now, $\alpha=\tau\circ((\sigma\ox\c N^*)\ox\Og^b_X)$ by definition, and
$\tau|Y$ is a perfect pairing because $X$ is smooth.  Hence, by the
preceding paragraph, $Y$ is invariant under $\eta$ if and only if there
exists a $\zeta$ such that $\bigwedge^a d\circ\zeta$ is equal to
$(\sigma\ox\c N^*)|Y$.  By the second paragraph, such a $\zeta$ exists
if and only if $Y$ is a solution of $\sigma=0$.
 \end{proof}

 \begin{proposition} Let $X$ be a projective scheme, $Y\subseteq X$ a
reduced closed subscheme of dimension $b$ such that the induced map
$\H^b(\Og^b_X)\to\H^b(\Og^b_Y)$ is nonzero.  Let $\eta\:\Og^b_X\to \c L$
be a Pfaff field, $S$ its singular locus.  Assume that no
$b$-dimensional component of $Y$ lies in $S$ and that $Y$ is invariant
under $\eta$.  Then $\H^b(\c I_{Y\cap S,\,X}\otimes\c L)\neq 0$.

Furthermore, if $\H^b(\eta)=0$, then the following three statements
hold:
\par\smallskip
  {\rm(1)}\enspace The restriction map $\H^{b-1}(\c L)\to\H^{b-1}(\c
  L|(Y\cap S))$ is not surjective.\par
  {\rm(2)}\enspace The intersection $Y\cap S$ has maximal dimension, 
$\dim (Y\cap S)=b-1$.\par
  {\rm(3)}\enspace The induced map $\H^b(\Og^b_Y)\to\H^b(\c L|Y)$ is
  surjective, but not bijective.
 \end{proposition}

 \begin{proof} Let $\mu\:\Og^b_Y\to\c L|Y$ be the map making Diagram
(3.1.2) 
commute. First, let us see that Diagram (3.1.2) induces a
commutative diagram
 $$
 \begin{CD}
 \Og^b_X @>\eta'>> \c I_{Y\cap S,\,X}\ox\c L\\
 @VVV @VVV\\
 \Og^b_Y @>\mu'>> \c I_{Y\cap S,\,Y}\ox\c L|Y
 \end{CD}\eqno(3.3.1)
 $$
 in which $\mu'$ is surjective.  Indeed, the image of $\eta$ is $\c
I_{S,\,X}\ox\c L$ owing to the definition of $S$.  So the image of
$\Og^b_X$ in $\c L|Y$ is $\c I_{Y\cap S,\,Y}\ox\c L|Y$.  Since $\Og^b_X
\to \Og^b_Y$ is surjective, the image of $\mu$ is $\c I_{Y\cap
S,\,Y}\ox\c L|Y$ too; whence, $\mu$ induces $\mu'$.  Finally, the
natural map $\c I_{S,\,X}\to \c I_{Y\cap S,\,Y}$ factors through $\c
I_{Y\cap S,\,X}$; whence, $\eta$ induces $\eta'$.

Since $Y$ is of dimension $b$ and reduced, $\Og^b_Y$ is invertible in
codimension 0.  Since no $b$-dimensional component of $Y$ lies in $S$,
also $\c I_{Y\cap S,\,Y}\ox\c L|Y$ is invertible in codimension 0.
Hence, since $\mu'$ is surjective, $\mu'$ is bijective in codimension 0.
So the support of its kernel has dimension at most $b-1$.  Therefore,
$\H^b(\mu')$ is bijective.  Now, $\H^b(\Og^b_X)\to\H^b(\Og^b_Y)$ is
nonzero.  Hence, since Diagram (3.3.1) is commutative,
        $$\H^b(\eta')\neq 0.\eqno(3.3.2)$$
 Thus $\H^b(\c I_{Y\cap S,\,X}\ox\c L)\neq 0$, as asserted.

Assume $\H^b(\eta)=0$ now.  Form the standard exact sequence
        $$0\to\c  I_{Y\cap S,\,X}\to\c O_X\to\c O_{Y\cap S}\to 0,$$
 tensor it with $\c L$, and extract the following exact sequence
of cohomology:
        $$\H^{b-1}(\c L)\to\H^{b-1}(\c L|(Y\cap S))
        \to\H^b(\c I_{Y\cap S,\,X}\ox\c L)\to\H^b(\c L).$$
 By exactness at $\H^b(\c I_{Y\cap S,\,X}\ox\c L)$, the image of
$\H^{b-1}(\c L|(Y\cap S))$ contains the image of $\H^b(\eta')$ since
$\H^b(\eta)=0$.  But, this image is nonzero owing to (3.3.2).  Hence, by
exactness at $\H^{b-1}(\c L|(Y\cap S))$, the first map is not
surjective; that is, (1) holds.

In particular, $\H^{b-1}(\c L|Y\cap S)\neq 0$.  Hence $\dim (Y\cap
S)\geq b-1$.  But no $b$-dimensional component of $Y$ lies in $S$.
Therefore, $\dim (Y\cap S)=b-1$; that is, (2) holds.

Since the cokernel of $\mu$ is supported on $Y\cap S$, the map
$\H^b(\mu)$ is surjective.  But it is not bijective since Diagram
(3.1.2) is commutative, since $\H^b(\eta)=0$ and since 
$\H^b(\Og^b_X)\to \H^b(\Og^b_Y)$ is nonzero.  Thus (3) holds.
 \end{proof}

\begin{proposition}
 Let $X$ be a smooth projective scheme in characteristic $0$. Let
$Y\subseteq X$ be a closed subscheme of pure dimension $b$, and $\c N$
its normal sheaf.  Assume ${\rm Pic}(X)= {\bf Z}$, and assume there is a
closed curve $C\subseteq X$ contained in the smooth locus $U$ of $Y$
such that $\deg\c N|C>0$.  Let $\eta\:\Og^b_X\to \c L$ be a Pfaff field,
$S$ its singular locus.  Assume no $b$-dimensional component of $Y$ lies
in $S$, and assume $Y$ is invariant under $\eta$.  Then $\dim (Y\cap
S)=b-1$.  \end{proposition}

\begin{proof} Since $Y$ is invariant, $\eta$ induces a map 
$\mu\:\Og^b_Y\to\c L|Y$ making (3.1.2) commute.  And $\mu$ is surjective
off $Z:=Y\cap S$.

Let $\c I_{Y,X}$ be the ideal of $Y$ in $X$, and consider the second
fundamental sequence,
 $$0\to\c I_{Y,X}/\c I^2_{Y,X}\to\Og^1_X|Y\to\Og^1_Y\to 0.\eqno(3.4.1)$$
 It is right exact, and is exact on the smooth locus $U$ of $Y$.

Set $n:=\dim X$ and $a:=n-b$. Then (3.4.1) induces an isomorphism
        $$\Og^b_U\risom\Og^n_X|U\ox\bigwedge^a\c N|U.$$ Set $V:=U-U\cap
Z$.  Then $\mu|V$ is a surjection between invertible sheaves, whence a 
bijection.  So there is an isomorphism
$\rho\:\Og^n_X|V\ox\bigwedge^a\c N|V\risom\c L|V$.

 Proceeding by way of contradiction, assume $\dim Z\le b-2$.  Then $U\cap
Z$ has codimension at least 2 in $U$.  Since $U$ is smooth and
$\Og^n_X|U\ox\bigwedge^a\c N|U$ and $\c L|U$ are invertible, $\rho$
extends to an isomorphism $\wt\rho\:\Og^n_X|U\ox\bigwedge^a\c
N|U\risom\c L|U$.

By hypothesis, $C\subseteq U$ and $\deg\c N|C>0$. So 
$\deg\c L|C>\deg\Og^n_X|C$. 
Since ${\rm Pic}(X)={\bf Z}$, the sheaf
$\c L\ox(\Og^n_X)^{-1}$ is ample. So, by the Kodaira vanishing theorem, 
$\h^b(\c L)=0$. Now, since the characteristic is 0, the natural map 
$\H^b(\Og^b_X)\to\H^b(\Og^b_Y)$ is nonzero by Proposition 2.1. 
Hence $\dim (Y\cap S)= b-1$ by Proposition 3.3,
a contradiction.
\end{proof}

\begin{lemma} Let $X$ be a projective scheme, $Y\subset X$ a reduced 
closed subscheme of pure dimension $b$.  Assume $Y$ is Gorenstein, and
has normal-crossings in codimension $1$.  Let $r$ be the number of
irreducible components of $Y$, and $\bog_Y$ its dualizing sheaf.  Let
$\eta\:\Og^b_X\to\c L$ be a Pfaff field, $S$ its singular locus.  Assume
no irreducible component of $Y$ lies in $S$, and assume $Y$ is invariant
under $\eta$.  Let $m\geq 0$.  Then $\h^0(\bog_Y^m\ox\c L^{-m}|Y)\leq
r$.
  \end{lemma}

\begin{proof} Since $Y$ is invariant, $\eta$ induces a map
$\mu\:\Og^b_Y\to\c L|Y$, which is surjective away from $Y\cap S$.  Now,
no irreducible component of $Y$ lies in $S$, and $Y$ is reduced of pure
dimension $b$.  Hence $\mu$ is generically bijective.

  Let $f\:\wt Y\to Y$ be the normalization map, $\lambda\:\Og^b_Y\to
f_*\Og^b_{\wt Y}$ the induced map.  Since $Y$ is reduced, whence
generically smooth, $\lambda$ is generically bijective.  Now, $Y$ has
normal crossings in codimension 1; hence, $\lambda$ is also surjective
in codimension 1.

 Let $\wt\gamma\:\Og^b_{\wt Y}\to\bog_{\wt Y}$ be the ``class'' map. It
is bijective on the smooth locus of $\wt Y$, so in codimension 1 as $\wt
Y$ is normal.  Hence the composition $\beta:=f_*\wt\gamma\circ\lambda$ is
generically bijective and surjective in codimension 1.  
Since the kernel of $\beta$ is torsion, but $\c L|Y$ is torsion free, 
it follows that $\mu$ factors through $\beta$ in codimension 1. 
Now, $Y$ is Gorenstein, so Cohen--Macaulay, and 
$\c L|Y$ is invertible.  Hence 
$\mu$ factors through $\beta$, yielding a map
$\tau\:f_*\bog_{\wt Y}\to\c L|Y$.  And $\tau$ is generically bijective.

Since $f$ is finite, $f_*\bog_{\wt Y}=\IHom(f_*\c O_{\wt Y},\bog_Y)$.
Now, $\bog_Y$ is dualizing.  Hence 
$f_*\c O_{\wt Y}=\IHom(f_*\bog_{\wt Y},\bog_Y)$. Applying 
$\IHom(\bullet,\bog_Y)$ to
$\tau$, we get a generically bijective map $\rho\:\IHom(\c L|Y,\bog_Y)\to
f_*\c O_{\wt Y}$.

Set $\c M:=\IHom(\c L|Y,\bog_Y)$, and consider the composition
 $$ \begin{CD} \iota\: \c M^m @>\rho^{\otimes m}>> (f_*\c
        O_{\wt Y})^{\otimes m} @>\pi>> f_*\c O_{\wt Y} \end{CD} $$
 where $\pi$ is given by multiplication.  Since $f$ is birational, $\pi$
 is generically bijective.  So, as $\rho$ is generically bijective,
 $\iota$ is too.  Hence, since $\c M^m$ is invertible and $Y$ is
 reduced, $\iota$ is globally injective.  Thus $\h^0(\c
 M^m)\leq\h^0(f_*\c O_{\wt Y})$.  But $\h^0(f_*\c O_{\wt Y})=\h^0(\c
 O_{\wt Y})=r$.  The assertion follows.  \end{proof}

\begin{theorem} Let $X$ be a smooth projective scheme with ${\rm
Pic}(X)={\bf Z}$.  Let $Y\subset X$ be a reduced closed subscheme, and
assume $Y$ is a Cartier divisor with normal crossings in
codimension~$1$.  Let $\sigma\in\H^0(\Og^1_X\ox\c N)$ be a Pfaff
equation, $S$ its singular locus.  Assume no irreducible component of
$Y$ lies in $S$, and assume $Y$ is a solution of $\sigma=0$.  Then $\deg
Y\le\deg\c N$; furthermore, $\deg Y<\deg\c N$ if\/ $Y$ is smooth in
codimension $1$ and the characteristic is $0$.
  \end{theorem}

\begin{proof} Since $X$ is smooth and $Y\subseteq X$ is a Cartier
divisor, $Y$ is Gorenstein, and its dualizing sheaf $\bog_Y$ is given by
the formula $\bog_Y=\Og^n_X(Y)|Y$ where $n:=\dim X$.

Let $\eta\:\Og^{n-1}_X\to\c L$ be the Pfaff field corresponding to
$\sigma$, so $\c L:=\Og^n_X\ox\c N$. Then
       $$ \bog_Y\ox\c L^{-1}|Y=\c N^{-1}(Y)|Y.\eqno(3.6.1) $$
 By Lemma 3.5, there is an integer $r$ such that $\h^0(\c
N^{-m}(mY)|Y)\leq r$ for every $m\geq 0$.  Hence $\c N^{-1}(Y)$ is
nonpositive because $\dim Y>0$. So $\deg Y\le\deg\c N$, as asserted.

Furthermore, since $Y$ is invariant, $\eta$ induces a map
$\mu\:\Og^{n-1}_Y\to\c L|Y$.  And $\mu$ is surjective off $Y\cap S$.  Now,
$\Og^{n-1}_Y$ is generically invertible because $Y$ is reduced, and
$Y-Y\cap S$ is dense in $Y$; hence, $\mu$ is generically injective.

Assume now that $Y$ is smooth in codimension 1.  Then the ``class'' map
$\gamma\:\Og^{n-1}_Y\to\bog_Y$ is an isomorphism in codimension 1.  On
the other hand, since $Y$ is invariant, $\eta$ induces a map
$\mu\:\Og^{n-1}_Y\to\c L|Y$.  Since $Y$ is Gorenstein, so
Cohen--Macaulay, and since $\c L|Y$ is invertible, $\mu$ factors through
$\gamma$, yielding a map $\tau\:\bog_Y\to\c L|Y$.  As $\mu$ is
generically injective, so is $\tau$.

Assume $\deg Y=\deg\c N$.  Then $\c N^{-1}(Y)=\c O_X$ since ${\rm
Pic}(X)={\bf Z}$.  So $\bog_Y=\c L|Y$ by (3.6.1).  Hence $\tau$
corresponds to a generically nonzero, everywhere regular function $f$ on
$Y$.  Since $Y$ is projective, $f$ is locally constant, so everywhere
nonzero. Thus $\tau$ is an isomorphism. Hence $\mu$ is, like $\gamma$,
an isomorphism in codimension 1.  Thus $\codim(Y\cap S,Y)\ge 2$.

However, $Y$ is projective and smooth in codimension 1; so the smooth
locus of Y contains a (smooth) projective curve by Bertini's
Theorem. Hence, by Proposition 3.4, the characteristic must be nonzero.
Thus the second assertion is proved.
\end{proof}

\section{Projective spaces}

\begin{definition} Let $X:=\IP^{n_1}\x\dots\x\IP^{n_s}$ and let
$\eta\:\Og^b_X\to\c L$ be a Pfaff field.  Call $\eta$ {\it fibered}
if the $\IP^{n_i}$ can be grouped to give a decomposition $X=X_1\x X_2$
such that
 \par\smallskip
  (4.1.1)\enspace $\dim X_1=b$, \par
  (4.1.2)\enspace $\c L=\Og^b_{X_1}\bx\c M_2$ where $\c M_2$ is nonnegative on
$X_2$, and \par
  (4.1.3)\enspace $\eta\:\Og^b_X\to\Og^b_{X_1}\bx\c O_{X_2}\to\c L$
 where the first map is the natural
 \indent\hphantom{(4.1.3)}\enspace surjection and the second arises from
a section of $\c M_2$.
 \end{definition}

\begin{lemma} Let $X:=\IP^{n_1}\x\dots\x\IP^{n_s}$ and let
$\eta\:\Og^b_X\to\c L$ be a Pfaff field.  Then $\eta$ is fibered if
and only if $\H^b(\c L)\neq 0$.
 \end{lemma}
 \begin{proof} Say $\c L=\c O_X(m_1,\dots,m_s)$.  
By the K\"unneth formula,
 $$
 \H^b(\c O_X(m_1,\dots,m_s))=
 \sum\H^{b_1}(\c O_{\IP^{n_1}}(m_1))\ox\cdots\ox\H^{b_s}(\c
 O_{\IP^{n_s}}(m_s))
 $$
 where the sum ranges over all $s$-tuples $(b_1,\dots,b_s)$ of
nonnegative integers $b_i$ such that $b_1+\cdots+b_s=b$.  So $\H^b(\c
O_X(m_1,\dots,m_s))\neq 0$ if and only if there is such an $s$-tuple
$(b_1,\dots,b_s)$ such that $\H^{b_i}(\c O_{\IP^{n_i}}(m_i))\neq 0$ for
all $i$.  However, $\H^{b_i}(\c O_{\IP^{n_i}}(m_i))\neq 0$ if and only
if either $b_i=n_i$ and $m_i\le-n_i-1$ or $b_i=0$ and $m_i\ge 0$.

Reorder the $\IP^{n_i}$ so that $m_i<0$ for $1\le i\le t$ and $m_i\ge 0$
for $t+1\le i\le s$.  Then $\H^b(\c L)\neq 0$ if and only if
        $$b= n_1+\cdots+n_t \text{  and }m_i\le-n_i-1
        \text{ for } 1\le i\le t.\eqno(4.2.1)$$
 (Up to this point, $\eta$ has played no role.)

Suppose $\eta$ is fibered.  Then, owing to (4.1.2), the $\IP^{n_i}$ must
be grouped into the first $t$ for $X_1$ and the rest for $X_2$.  Also,
(4.1.1) and (4.1.2) yield (4.2.1).  Hence $\H^b(\c L)\neq 0$.

Conversely, suppose $\H^b(\c L)\neq 0$; then (4.2.1) holds.  Group the
$\IP^{n_i}$ into the first $t$ for $X_1$ and the rest for $X_2$ 
to get $X=X_1\x X_2$.  Then (4.1.1) holds.

According to Definition 3.1, the field $\eta$ corresponds to a nonzero
section
 $$
 \sigma\in \H^0\bigl(\Og^a_X(m_1+n_1+1,\dots,m_s+n_s+1)\bigr)
 \hbox{ where } a:=n_1+\cdots+n_s-b.  
 $$
 Now,  $\Og^1_X$ is equal to the sum of the
pullbacks of the sheaves $\Og^1_{\IP^{n_i}}$.  So 
 $$
 \Og^{a}_X=\sum\Og^{e_1}_{\IP^{n_1}}\bx\cdots\bx\Og^{e_s}_{\IP^{n_s}}
 $$
 where the sum ranges over all $s$-tuples $(e_1,\dots,e_s)$ of
nonnegative integers $e_i$ such that $e_1+\cdots+e_s=a$.  Via the
K\"unneth formula, $\sigma$ becomes an element of the sum
 $$
 \sum\H^0\bigl(\Og^{e_1}_{\IP^{n_1}}(m_1+n_1+1)\bigr)\ox\cdots\ox
 \H^0\bigl(\Og^{e_s}_{\IP^{n_s}}(m_s+n_s+1)\bigr).\eqno(4.2.2)
 $$

Since $\sigma\neq 0$, there is an $s$-tuple $(e_1,\dots,e_s)$ such that
 $$
 \H^0\bigl(\Og^{e_i}_{\IP^{n_i}}(m_i+n_i+1)\bigr)\neq 0\text{ for all }i.
 $$
This $s$-tuple is unique; in fact, let's now show that
 \par\smallskip
 (4.2.3)\enspace $e_i=0$ and $m_i=-n_i-1$ for $1\le i\le t$, and\par
 (4.2.4)\enspace $e_i=n_i$ for $t+1\le i\le s$.
 \par\smallskip
 Indeed, first fix $i\le t$.  Then $m_i+n_i+1\leq 0$ by (4.2.1).  Now,
$\Og^{e_i}_{\IP^{n_i}}$ embeds in a direct sum of copies of $\c
O_{\IP^{n_i}}(-e_i)$; hence,
$\H^0\bigl(\Og^{e_i}_{\IP^{n_i}}(m_i+n_i+1)\bigr)$ is nonzero only if
$m_i+n_i+1- e_i\ge0$.  But, $e_i\ge0$.  Therefore, (4.2.3) holds.

Owing to  (4.2.3), we have $e_{t+1}+\cdots+e_s=a$.  However,
$a=n_{t+1}+\cdots+n_s$ again owing to (4.2.1).  Since $e_i\le n_i$ for
all $i$, therefore (4.2.4) holds.

By (4.2.1) and (4.2.3), $\c L=\Og^b_{X_1}\bx\c M_2$, where 
$\c M_2:=\c O_{X_2}(m_{t+1},\dots,m_s)$.  Since $m_i\ge 0$ for 
$t+1\le i\le s$, the sheaf $\c M_2$ is nonnegative.  So (4.1.2) holds

Owing to (4.2.3) and (4.2.4), the sum in (4.2.2) reduces to the single
term
 \begin{multline*}
 \H^0(\c O_{\IP^{n_1}})\ox\cdots\ox\H^0(\c O_{\IP^{n_t}})\\
 \ox\H^0\bigl(\Og^{n_{t+1}}_{\IP^{n_{t+1}}}(m_{t+1}+n_{t+1}+1)\bigr)
  \ox\cdots\ox\H^0\bigl(\Og^{n_s}_{\IP^{n_s}}(m_s+n_s+1)\bigr).
 \end{multline*}
 This term is, by the   K\"unneth formula, equal to
        $$\H^0\bigl(\c O_{X_1}\bx\Og^a_{X_2}(m_{t+1}+n_{t+1}+1,\ldots,
        m_s+n_s+1)\bigr). $$
Since $a=n_{t+1}+\cdots+n_s$, the latter corresponds to $\H^0(\c M_2)$.
Thus 
$\sigma$ corresponds to a section of $\c M_2$, and (4.1.3) follows.
\end{proof}

 \begin{theorem} Set $X:=\IP^{n_1}\x\dots\x\IP^{n_s}$.  Let $Y\subseteq
X$ be a reduced closed subscheme of dimension $b$ and of multidegree
$(d_1,\ldots,d_s)$.  Assume either the characteristic is $0$, or it is
$p>0$ and $p$ does not divide some $d_i$.  Let $\eta\: \Og^b_X\to\c L$
be a Pfaff field, and  $S$ its singular locus.  Assume  no
$b$-dimensional component of $Y$ lies in $S$, and $Y$ is invariant
under $\eta$.  \par\smallskip
  {\rm(1)}\enspace Then $\H^b(\c I_{Y\cap S,\,X}\otimes\c L)\neq 0$ where $\c
I_{Y\cap S,\,X}$ is the ideal of $Y\cap S$ in $X$.\par
  {\rm(2)}\enspace If $\eta$ is nonfibered, then $\dim (Y\cap S)=b-1$.\par
  {\rm(3)}\enspace If $\eta$ is nonfibered and $\h^b(\Og^b_Y)=1$, then
  $\h^b(\c L|Y)=0$.
 \end{theorem}
 \begin{proof}
 The natural map $\H^b(\Og^b_X)\to\H^b(\Og^b_Y)$ is nonzero either by
Proposition~2.1 or by Proposition 2.2.  Hence Proposition 3.3 yields
(1).

Assume $\eta$ is nonfibered.  Then $\H^b(\c L)= 0$ by Lemma 4.2,
so $\H^b(\eta)=0$.  Hence Proposition 3.3 yields (2).  It also yields
that $\h^b(\Og^b_Y)>\h^b(\c L|Y)$; whence, (3) holds.
 \end{proof}

\begin{remark} Let $X:=\IP^{n_1}\x\dots\x\IP^{n_s}$, and 
$Y\subseteq X$ be a reduced connected closed subscheme of pure dimension
$b$.
Let $\eta\:\Og^b_X\to\c L$ be a fibered Pfaff field, and $S$ its singular 
locus.  Note that (4.1.3) implies $S=X_1\x S_2$ where $S_2$ is a divisor
of
$X_2$.  Assume no irreducible component of $Y$ lies in $S$, and $Y$
is invariant under $\eta$.

Suppose the characteristic is 0.  Let $Z$ be the image of $Y$ on $X_2$ 
under the projection map, and $\pi\:Y\to Z$ be the induced surjection. 
Since $Y$ is invariant under $\eta$, by (4.1.3) the differential $d\pi_y$ 
is zero 
at $y\in Y-Y\cap S$.  Now, since the characteristic is 0, the induced map 
$\pi\:Y\to Z$ is generically smooth.  As $Y-Y\cap S$ is dense in $Y$, 
the map $d\pi_y$ is zero for a general $y\in Y$.  Hence $\dim Z=0$, and 
since $Y$ is connected, $Z=\{P\}$ for $P\in X_2-S_2$. 
So $Y=X_1\x\{P\}$.  Hence $Y\cap S$ is empty; so (2) of
Theorem 4.3 fails.  Furthermore, $\c L|Y=\Og^b_Y$; so (3) fails as well. 
 \end{remark}

 \begin{corollary} Let $X:=\IP^{n}$, and $Y\subseteq X$ a reduced
closed subscheme of dimension $b$ and of degree $d$.  Assume either the
characteristic is $0$ or it is $p>0$ and $p$ does not divide $d$.  Let
$\eta\: \Og^b_X\to\c L$ be a Pfaff field, and $S$ its singular locus.
Assume no $b$-dimensional component of $Y$ lies in $S$, and $Y$ is
invariant under $\eta$.  Set $l:=\deg\c L$.  Then
     $$ \reg(Y\cap S)>  b+l \text{ and }\dim (Y\cap S)=b-1$$
 where $\reg(Y\cap S)$ denotes the Castelnuovo--Mumford regularity.

Assume $\h^b(\Og^b_Y)=1$ in addition.  Then $\h^b(\c L|Y)=0$.
Furthermore, if $Y$ is arithmetically Cohen--Macaulay, then $\reg Y\leq
b+l+1$; if $Y$ is $s$-subcanonical, then $s\le l-1$.
 \end{corollary}

 \begin{proof} Note that $\eta$ is nonfibered by Definition 4.1 since
$0<b<n$.  Hence the first three assertions are immediate
consequences of those of Theorem~4.3.

Set $m:=b+l+1$.  Then $m\ge0$.  Indeed, $\eta$ corresponds to a nonzero
section of $\Og_X^{n-b}(l+n+1)$. Since $\Og_X^{n-b}$ embeds in a direct
sum of copies of $\c O_X(b-n)$, we get a nonzero section of 
$\c O_X(m)$. So $m\ge0$.

By definition, $\reg Y\leq m$ if $\H^i(\c I_{Y,\,X}(m-i))$ vanishes for
$i\ge 1$.  

Consider the standard exact sequence
        $$\H^{i-1}(\c O_X(m-i))\to\H^{i-1}(\c O_Y(m-i))
        \to\H^i(\c I_{Y,\,X}(m-i))\to\H^i(\c O_X(m-i)).\eqno(4.5.1)$$
 Since $m\ge0$, the last group vanishes for $i\ge 1$.  For 
$i=b+1$, the second group is $\H^b(\c L|Y)$, hence zero by the third 
assertion. So $\H^i(\c I_{Y,\,X}(m-i))$ vanishes for $i=b+1$.

Suppose $Y$ is arithmetically Cohen--Macaulay.  Then the first map 
in (4.5.1) is surjective for $i=1$, and $\H^{i-1}(\c O_Y(m-i))$ vanishes 
for $i\neq 1,\,b+1$. So $\H^i(\c I_{Y,\,X}(m-i))$ vanishes for 
$i\neq 0,\,b+1$ as well. Thus $\reg Y\leq m$.

Finally, assume $Y$ is $s$-subcanonical; that is, $\c O_Y(s)$ is its
dualizing sheaf.  Now, $\h^b(\c L|Y)=0$, and $\c L|Y=\c O_Y(l)$.  Hence
$\h^0(\c O_Y(s-l))=0$ by duality.  So $s-l\le-1$.
 \end{proof}

 \begin{remark} Let $X:=\IP^{n}$, and $Y\subseteq X$ a reduced closed
subscheme of dimension $b$.  Assume $\h^b(\Og^b_Y)=1$.  Let  $\bog_Y$ be
the dualizing sheaf.  
Then, by duality,
        $$\H^b(\Og^b_Y)=\Hom(\Og^b_Y,\bog_Y)^* \text{ and }
        \H^b(\Og^b_Y(1))=\Hom(\Og^b_Y,\bog_Y(-1))^*.$$
  Take any nonzero map $\gamma\:\Og^b_Y\to\bog_Y$, and let $T$ be the
support of its cokernel.  Take a hyperplane $H$ that doesn't contain
$T$.  Then the image of $\gamma$ is not contained in $\bog_Y(-H)$.  So
$\gamma\notin \Hom(\Og^b_Y,\bog_Y(-H))$.  Hence
$\Hom(\Og^b_Y,\bog_Y(-1))=0$.  Thus $\h^b(\Og^b_Y(1))=0$.

Let $\eta\: \Og^b_X\to\c L$ be a Pfaff field, and $S$ its singular
locus.  Assume $Y$ is invariant under $\eta$, and form the map
$\mu\:\Og^b_Y\to\c L|Y$ induced by $\eta;$ see Diagram (3.1.2).  Assume
no $b$-dimensional component of $Y$ lies in $S$.  Then $\mu$ is
surjective off $Y\cap S$; furthermore, $\dim Y\cap S\le b-1$.  Since
$\h^b(\Og^b_Y(1))=0$, it follows that $\h^b(\c L(1)|Y)=0$.

Set $l:=\deg\c L$.  The proof of Corollary 4.5 now yields these two
weaker conclusions: if $Y$ is arithmetically Cohen-Macaulay, then $\reg
Y\leq b+l+2$; if $Y$ is $s$-subcanonical, then $s\leq l$.

For the stronger conclusions of Corollary 4.5 to hold however, it is
necessary that either the characteristic be $0$ or  it be $p>0$ and
$p$ not divide $d$ where $d:=\deg Y$.  For example, take $n:=2$ and
$p>0$, and let $Y$ be a smooth curve such that $p$ divides $d$.  Then it
is possible to find a Pfaff field $\eta\: \Og^1_X\to\c L$ where $\c L :=
\c O_X(d-3)$ such that $Y$ is invariant under 
$\eta$;
see Remark 14 on p.~8 of \cite{E}.
 \end{remark}

 \begin{remark} Let $Y$ be an integral projective scheme of
dimension $b$.  Suppose that $Y$ has
normal-crossings singularities in codimension 1; that is, there is a
closed set $F$ of dimension at most $b-2$ such that $Y-F$ has two sheets
meeting transversally.  Let's show $\h^b(\Og^b_Y)=1$.

Indeed, let $f\:Y^*\to Y$ be the normalization map.  The natural map
$\Og_Y^1\to f_*\Og_{Y^*}^1$ is bijective off $F$.  Since $\dim F\le
b-2$, the induced map $\H^b(\Og^b_Y) \to \H^b(f_*\Og_{Y^*}^b)$ is
bijective.  Its target is equal to $\H^b(\Og_{Y^*}^b)$ as $f$ is
finite.  Thus we may assume $Y$ is normal.

Consider the ``class'' map from $\Og_Y^b$ to the dualizing sheaf
$\bog_Y$.  This map is bijective on the smooth locus of $Y$.  Since the
singular locus of $Y$ has dimension at most $b-2$, because $Y$ is
normal, the class map induces a bijection from $\H^b(\Og^b_Y)$ to
$\H^b(\bog_Y)$.  The latter group is, by duality, equal to $H^0(\c
O_Y)$.  Since $Y$ is integral, $h^0(\c O_Y)=1$.  Thus $\h^b(\Og^b_Y)=1$.
 \end{remark}

\section*{Acknowledgments}

The first author thanks A.\ Campillo, L.\ G.\ Mendes, P.\ Sad, 
M.\ Soares, and especially J.\ V.\ Pereira for helpful discussions on the 
subject.  
He is also grateful to CNPq for a grant, Proc.\ 202151/90-5, supporting
a year-long visit to MIT, and grateful to MIT for its hospitality. 
He was also supported by PRONEX, Conv\^enio 41/96/0883/00, CNPq, 
Proc.\ 300004/95-8, and FAPERJ, Proc.{} E-26/170.418/2000-APQ1.

 The second author thanks IMPA, Rio de Janeiro, ICMC-USP, S\~ao Carlos,
and the XVI and XVII Escolas de \'Algebra, Bras\'\i lia and Cabo Frio,
for their invitations and financial support, which enabled this work
to be initiated, pursued, and presented.

\end{document}